\documentclass{dcds}
\usepackage{epsfig}
\usepackage{amsfonts,amsmath,amsthm}
\usepackage[latin1]{inputenc}
\def\currentvolume{}
 \def\currentissue{}
  \def\currentyear{2005}
   \def\currentmonth{}
    \def\ppages{1--??}

\newtheorem{theorem}{Theorem}[section]
\newtheorem{definition}{Definition}[section]
\newtheorem{proposition}[theorem]{Proposition}
\newtheorem{corollary}[theorem]{Corollary}
\newtheorem{lemma}[theorem]{Lemma}
\newtheorem{example}{Example}[section]
\newtheorem{remark}[theorem]{Remark}
\renewcommand{\theequation}{\thesection.\arabic{equation}}

\title[Domains of attraction for autonomous discrete dynamical
systems]{Methods for determination and approximation of domains of
attraction in the case of autonomous discrete dynamical systems}
\author[St. Balint, E. Kaslik, A.M. Balint, A. Grigis]{}
\email{balint@balint.math.uvt.ro}



\begin{document}

\maketitle

\centerline{\scshape  Stefan Balint}
\medskip
{\footnotesize \centerline{ Department of Mathematics }
\centerline{ West University of Timi\c{s}oara} \centerline{ Bd. V.
Parvan nr. 4, 300223, Timi\c{s}oara, Romania } }

\medskip

\centerline{\scshape  Eva Kaslik}
\medskip
{\footnotesize \centerline{ Department of Mathematics }
\centerline{ West University of Timi\c{s}oara} \centerline{ Bd. V.
Parvan nr. 4, 300223, Timi\c{s}oara, Romania } } {\footnotesize
\centerline{ and }\centerline{ L.A.G.A, UMR 7539, Institut
Galil\'{e}e } \centerline{ Universit\'{e} Paris 13} \centerline{
99 Avenue J.B. Cl\'{e}ment, 93430, Villetaneuse, France } }

\medskip

\centerline{\scshape  Agneta Maria Balint}
\medskip
{\footnotesize \centerline{ Department of Physics } \centerline{
West University of Timi\c{s}oara} \centerline{ Bd. V. Parvan nr.
4, 300223, Timi\c{s}oara, Romania } }

\medskip

\centerline{\scshape  Alain Grigis}
\medskip
{\footnotesize \centerline{ L.A.G.A, UMR 7539, Institut
Galil\'{e}e } \centerline{ Universit\'{e} Paris 13} \centerline{
99 Avenue J.B. Cl\'{e}ment, 93430, Villetaneuse, France } }

\medskip

\begin{abstract}
A method for determination and two methods for approximation of
the domain of attraction $D_{a}(0)$ of an asymptotically stable
steady state of an autonomous, $\mathbb{R}$-analytical, discrete
system is presented. The method of determination is based on the
construction of a Lyapunov function $V$, whose domain of
analyticity is $D_{a}(0)$. The first method of approximation uses
a sequence of Lyapunov functions $V_{p}$, which converges to the
Lyapunov function $V$ on $D_{a}(0)$. Each $V_{p}$ defines an
estimate $N_{p}$ of $D_{a}(0)$. For any $x\in D_{a}(0)$ there
exists an estimate $N_{p^{x}}$ which contains $x$. The second
method of approximation uses a ball $B(R)\subset D_{a}(0)$ which
generates the sequence of estimates $M_{p}=f^{-p}(B(R))$. For any
$x\in D_{a}(0)$ there exists an estimate $M_{p^{x}}$ which
contains $x$. The cases $\|\partial_{0}f\|<1$ and
$\rho(\partial_{0}f)<1$ are treated separately (even though the
second case includes the first one) because significant
differences occur.
\end{abstract}

\newpage
\section{Introduction}
Let be the following discrete dynamical system:
\begin{equation}
\label{dyn.sys}
x_{k+1}=f(x_{k})\qquad k=0,1,2...
\end{equation}
where $f:\Omega\rightarrow\Omega$ is an $\mathbb{R}$-analytic
function defined on a domain $\Omega\subset\mathbb{R}^{n}$,
$0\in\Omega$ and $f(0)=0$, i.e. $x=0$ is a steady state (fixed
point) of (\ref{dyn.sys}).

For $r>0$, denote by $B(r)=\{x\in\mathbb{R}^{n}:\|x\|<r\}$ the
ball of radius $r$.

The steady state $x=0$ of (\ref{dyn.sys}) is "stable" provided
that given any ball $B(\varepsilon)$, there is a ball $B(\delta)$
such that if $x\in B(\delta)$ then $f^{k}(x)\in B(\varepsilon)$,
for $k=0,1,2,...$ \cite{Kelley-Peterson}.

If in addition there is a ball $B(r)$ such that
$f^{k}(x)\rightarrow 0$ as $k\rightarrow\infty$ for all $x\in
B(r)$ then the steady state $x=0$ is "asymptotically stable"
\cite{Kelley-Peterson}.

The domain of attraction $D_{a}(0)$ of the asymptotically stable
steady state $x=0$ is the set of initial states $x\in \Omega$ from
which the system converges to the steady state itself i.e.
\begin{equation}\label{da}
    D_{a}(0)=\{x\in\Omega | f^{k}(x)\stackrel{k\rightarrow\infty}{\longrightarrow}0\}
\end{equation}

Theoretical research shows that the $D_{a}(0)$ and its boundary
are complicated sets. In most cases, they do not admit an explicit
elementary representation. The domain of attraction of an
asymptotically stable steady state of a discrete dynamical system
is not necessarily connected (which is the case for continuous
dynamical systems). This fact is shown by the following example.

\begin{example}
\emph{Let be the function $f:\mathbb{R}\rightarrow\mathbb{R}$
defined by
$f(x)=\frac{1}{2}x-\frac{1}{4}x^{2}+\frac{1}{2}x^{3}+\frac{1}{4}x^{4}$.
The domain of attraction of the asymptotically stable steady state
$x=0$ is $D_{a}(0)=(-2.79,-2.46)\cup (-1,1)$ which is not
connected.}\qed
\end{example}

Different procedures are used for the approximation of the
$D_{a}(0)$ with domains having a simpler shape. For example, in
the case of the Theorem 4.20 pg. 170 \cite{Kelley-Peterson} the
domain which approximates the $D_{a}(0)$ is defined by a Lyapunov
function $V$ built with the matrix $\partial_{0}f$ of the
linearized system in $0$, under the assumption
$\|\partial_{0}f\|<1$. In \cite{KBBB1}, a Lyapunov function $V$ is
presented in the case when the matrix $\partial_{0}f$ is a
contraction, i.e. $\|\partial_{0}f\|<1$. The Lyapunov function $V$
is built using the whole nonlinear system, not only the matrix
$\partial_{0}f$. $V$ is defined on the whole $D_{a}(0)$, and more,
the $D_{a}(0)$ is the natural domain of analyticity of $V$. In
\cite{KBGB5}, this result is extended for the more general case
when $\rho(\partial_{0}f)<1$ (where $\rho(\partial_{0}f)$ denotes
the spectral radius of $\partial_{0}f$). This last result is the
following:

\begin{theorem} (see \cite{KBGB5})
If the function $f$ satisfies the following conditions:
\begin{equation}
  f(0) = 0
\end{equation}
\begin{equation}
  \rho(\partial_{0}f) < 1
\end{equation}
then $0$ is an asymptotically stable steady state. $D_{a}(0)$ is
an open subset of $\Omega$ and coincides with the natural domain
of analyticity of the unique solution $V$ of the iterative first
order functional equation
\begin{equation}
\label{ecV}
\begin{array}{ll}
\left\{\begin{array}{l}
V(f(x))-V(x)=-\|x\|^{2}\\
V(0)=0
\end{array}\right.
\end{array}
\end{equation}
The function $V$ is positive on $D_{a}(0)$ and
$V(x)\stackrel{x\rightarrow x^{0}}{\longrightarrow}+\infty$, for
any $x^{0}\in \partial D_{a}(0)$ ($\partial D_{a}(0)$ denotes the
boundary of $D_{a}(0)$) or for $\|x\|\rightarrow \infty$.
\label{theorem.DA.raza.subun}

The function $V$ is given by
\begin{equation}\label{formula.V}
    V(x)=\sum\limits_{k=0}^{\infty}\|f^{k}(x)\|^{2}\qquad
    \textrm{for any }x\in D_{a}(0)
\end{equation}
\end{theorem}

The Lyapunov function $V$ can be found theoretically using
relation (\ref{formula.V}). In the followings, we will shortly
present the procedure of determination and approximation of the
domain of attraction using the function $V$ presented in
\cite{KBBB1,KBGB5}.

The region of convergence $D_{0}$ of the power series development
of $V$ in $0$ is a part of the domain of attraction $D_{a}(0)$. If
$D_{0}$ is strictly contained in $D_{a}(0)$, then there exists a
point $x^{0}\in\partial D_{0}$ such that the function $V$ is
bounded on a neighborhood of $x^{0}$. Let be the power series
development of $V$ in $x^{0}$. The domain of convergence $D_{1}$
of the series centered in $x^{0}$ gives a new part
$D_{1}\setminus(D_{0}\bigcap D_{1})$ of the domain of attraction
$D_{a}(0)$. At this step, the part $D_{0}\bigcup D_{1}$ of
$D_{a}(0)$ is obtained.

If there exists a point $x^{1}\in\partial (D_{0}\bigcup D_{1})$
such that the function $V$ is bounded on a neighborhood of
$x^{1}$, then the domain $D_{0}\bigcup D_{1}$ is strictly included
in the domain of attraction $D_{a}(0)$. In this case, the
procedure described above is repeated in $x^{1}$.

The procedure cannot be continued in the case when it is found
that on the boundary of the domain $D_{0}\bigcup D_{1}\bigcup
...\bigcup D_{p}$ obtained at step $p$, there are no points having
neighborhoods on which $V$ is bounded.

This procedure gives an open connected estimate $D$ of the domain
of attraction $D_{a}(0)$. Note that $f^{-k}(D)$, $k\in\mathbb{N}$
is also an estimate of $D_{a}(0)$, which is not necessarily
connected.

The procedure described above is illustrated by the following
examples:

\begin{example}
\emph{Let be the $f:\mathbb{R}\rightarrow\mathbb{R}$ defined by
$f(x)=x^{2}$. Due to the equality $f^{k}(x)=x^{2^{k}}$ the domain
of attraction of the asymptotically stable steady state $x=0$ is
$D_{a}(0)=(-1,1)$. The Lyapunov function is
$V(x)=\sum\limits_{k=0}^{\infty}x^{2^{k+1}}$. The domain of
convergence of the series is $D_{0}=(-1,1)$ which coincides with
$D_{a}(0)$.}\qed
\end{example}

\begin{example}
\emph{Let be the function $f:\Omega=(-\infty,1)\rightarrow\Omega$
defined by $f(x)=\frac{x}{e+(1-e)x}$. Due to the equality
$f^{k}(x)=\frac{x}{e^{k}+(1-e^{k})x}$ the domain of attraction of
the asymptotically stable steady state $x=0$ is
$D_{a}(0)=(-\infty,1)$. The power series expansion of the Lyapunov
function $V(x)=\sum\limits_{k=0}^{\infty}|f^{k}(x)|^{2}$ in $0$ is
\begin{equation}\label{ex2.serie1}
    V(x)=\sum\limits_{m=2}^{\infty}(m-1)\sum\limits_{k=0}^{\infty}e^{-2k}(1-e^{-k})^{m-2}x^{m}
\end{equation}
The radius of convergence of the series (\ref{ex2.serie1}) is
\begin{equation}
r_{0}=\lim\limits_{m\rightarrow\infty}\sqrt[m]{(m-1)\sum\limits_{k=0}^{\infty}e^{-2k}(1-e^{-k})^{m-2}}=1
\end{equation}
therefore the domain of convergence of the series
(\ref{ex2.serie1}) is $D_{0}=(-1,1)\subset D_{a}(0)$. More,
$V(x)\rightarrow \infty$ as $x\rightarrow 1$ and $V(-1)<\infty$.
The radius of convergence of the power series expansion of $V$ in
$-1$ is
\begin{equation}
r_{-1}=\lim\limits_{m\rightarrow\infty}\sqrt[m]{\sum\limits_{k=1}^{\infty}\frac{e^{k}(e^{k}-1)^{m-2}[(m-3)e^{k}+2]}
{(2e^{k}-1)^{m+2}}}=1
\end{equation}
therefore, the domain of convergence of the power series
development of $V$ in $-1$ is $D_{-1}=(-2,0)$ which gives a new
part of $D_{a}(0)$.}\qed
\end{example}

Numerical results for more complex examples are given in
\cite{KBBB1,KBGB5}.

\section{Theoretical results when matrix $A=\partial_0 f$ is a contraction, i.e.
$\|A\|<1$}

The function $f$ can be written as
\begin{equation}\label{descomp.f}
    f(x)=Ax+g(x)\qquad\textrm{for any }x\in \Omega
\end{equation}
where $A=\partial_{0}f$ and $g:\Omega\rightarrow\Omega$ is an
$\mathbb{R}$-analytic function such that $g(0)=0$ and
$\lim\limits_{x\rightarrow 0}\frac{\|g(x)\|}{\|x\|}=0$.

\begin{proposition}
If $\|A\|<1$, then there exists $r>0$ such that
$B(r)\subset\Omega$ and $\|f(x)\|<\|x\|$ for any $x\in
B(r)\setminus\{0\}$. \label{proposition.ineg.norma}
\end{proposition}

\begin{proof}
Due to the fact that $\lim\limits_{x\rightarrow
0}\frac{\|g(x)\|}{\|x\|}=0$ there exists $r>0$ such that
$B(r)\subset\Omega$ and
\begin{equation}\label{B(R).definitie}
    \|g(x)\|<(1-\|A\|)\|x\|\qquad\textrm{for any } x\in B(r)\setminus\{0\}
\end{equation}
Let be $x\in B(r)\setminus\{0\}$. Inequality
(\ref{B(R).definitie}) provides that
\begin{equation}
    \|f(x)\|=\|Ax+g(x)\|\leq
    \|A\|\|x\|+\|g(x)\|<(\|A\|+1-\|A\|)\|x\|=\|x\|
\end{equation}
therefore, $\|f(x)\|<\|x\|$.
\end{proof}

\begin{definition}
Let be $R>0$ the largest number such that $B(R)\subset\Omega$ and
$\|f(x)\|<\|x\|$ for any $x\in B(R)\setminus\{0\}$.

If for any $r>0$ we have that $B(r)\subset\Omega$ and
$\|f(x)\|<\|x\|$ for any $x\in B(r)\setminus\{0\}$, then
$R=+\infty$ and $B(R)=\Omega=\mathbb{R}^{n}$.
\end{definition}

\begin{lemma}
\begin{itemize}
\item[a.] $B(R)$ is invariant to the flow of system
(\ref{dyn.sys}). \item[b.] For any $x\in B(R)$, the sequence
$(\|f^{k}(x)\|)_{k\in\mathbb{N}}$ is decreasing. \item[c.] For any
$p\geq 0$ and $x\in B(R)\setminus\{0\}$, $\Delta
V_{p}(x)=V_{p}(f(x))-V_{p}(x)<0$, where
\begin{equation}\label{formula.Vp}
    V_{p}(x)=\sum\limits_{k=0}^{p}\|f^{k}(x)\|^{2}\qquad
    \textrm{for }x\in \Omega
\end{equation}
\end{itemize}
\label{lem.B(R).proprietati}
\end{lemma}

\begin{proof}
\begin{itemize}
\item[a.] If $x=0$, then $f^{k}(0)=0$, for any $k\in\mathbb{N}$.
For $x\in B(R)\setminus\{0\}$, we have $\|f(x)\|<\|x\|$, which
implies that $f(x)\in B(R)$, i.e. $B(R)$ is invariant to the flow
of system (\ref{dyn.sys}).

\item[b.] By induction, it results that for $x\in B(R)$ we have
$f^{k}(x)\in B(R)$ and $\|f^{k+1}(x)\|\leq\|f^{k}(x)\|$, which
means that the sequence $(\|f^{k}(x)\|)_{k\in\mathbb{N}}$ is
decreasing.

\item[c.] In particular, for $p\geq 0$ and $x\in B(R)$, we have
$\|f^{p+1}(x)\|\leq\|f(x)\|<\|x\|$ and therefore, $\Delta
V_{p}(x)=\|f^{p+1}(x)\|^{2}-\|x\|^2<0$.
\end{itemize}
\end{proof}

\begin{corollary}
For any $p\geq 0$, there exists a maximal domain
$G_{p}\subset\Omega$ such that $0\in G_{p}$ and for $x\in
G_{p}\setminus\{0\}$, the (positive definite) function $V_{p}$
verifies $\Delta V_{p}(x)<0$. In other words, for any $p\geq 0$
the function $V_{p}$ defined by (\ref{formula.Vp}) is a Lyapunov
function for (\ref{dyn.sys}) on $G_{p}$. More, $B(R)\subset G_{p}$
for any $p\geq 0$.
\end{corollary}

\begin{theorem}
    $B(R)$ is an invariant set included in the domain of attraction $D_{a}(0)$.
\label{theorem.B(R).parte.DA}
\end{theorem}

\begin{proof}
Let be $x\in B(R)\setminus\{0\}$. We have to prove that
$\lim\limits_{k\rightarrow\infty}f^{k}(x)=0$.

The sequence $(f^{k}(x))_{k\in\mathbb{N}}$ is bounded: $f^{k}(x)$
belongs to $B(R)$. Let be $(f^{k_{j}}(x))_{j\in\mathbb{N}}$ a
convergent subsequence and let be
$\lim\limits_{j\rightarrow\infty}f^{k_{j}}(x)=y^{0}$. It is clear
that $y^{0}\in B(R)$.

It can be shown that
\begin{equation}\label{ineg0}
    \|f^{k}(x)\|\geq \|y^{0}\| \qquad\textrm{for any } k\in\mathbb{N}
\end{equation}
For this, observe first that $f^{k_{j}}(x)\rightarrow y^{0}$ and
$(\|f^{k_{j}}(x)\|)_{k\in\mathbb{N}}$ is decreasing (Lemma
\ref{lem.B(R).proprietati}). These imply that
$\|f^{k_{j}}(x)\|\geq \|y^{0}\|$ for any $k_{j}$. On the other
hand, for any $k\in\mathbb{N}$, there exists $k_{j}\in\mathbb{N}$
such that $k_{j}\geq k$. Therefore, as the sequence
$(\|f^{k}(x)\|)_{k\in\mathbb{N}}$ is decreasing (Lemma
\ref{lem.B(R).proprietati}), we obtain that
$\|f^{k}(x)\|\geq\|f^{k_{j}}(x)\| \geq \|y^{0}\|$.

We show now that $y^{0}=0$. Suppose the contrary, i.e. $y^{0}\neq
0$.

Inequality (\ref{ineg0}) becomes
\begin{equation}\label{contrad}
     \|f^{k}(x)\|\geq \|y^{0}\|>0 \qquad\textrm{for any } k\in\mathbb{N}
\end{equation}
By means of Lemma \ref{lem.B(R).proprietati}, we have that
$\|f(y^{0})\|<\|y^{0}\|$.

Therefore, there exists a neighborhood $U_{f(y^{0})}\subset B(R)$
of $f(y^{0})$ such that for any $z\in U_{f(y^{0})}$ we have
$\|z\|<\|y^{0}\|$. On the other hand, for the neighborhood
$U_{f(y^{0})}$ there exists a neighborhood $U_{y^{0}}\subset B(R)$
of $y^{0}$ such that for any $y\in U_{y^{0}}$, we have $f(y)\in
U_{f(y^{0})}$. Therefore:
\begin{equation}\label{ineg1}
    \|f(y)\|<\|y^{0}\|\qquad \textrm{for any } y\in U_{y^{0}}
\end{equation}
As $f^{k_{j}}(x)\rightarrow y^{0}$, there exists $\bar{j}$ such
that $f^{k_{j}}(x)\in U_{y^{0}}$, for any $j\geq \bar{j}$. Making
$y=f^{k_{j}}(x)$ in (\ref{ineg1}), it results that
\begin{equation}\label{ineg2}
    \|f^{k_{j}+1}(x)\|=\|f(f^{k_{j}}(x))\|<\|y^{0}\|
    \qquad \textrm{for }j\geq \bar{j}
\end{equation}
which contradicts (\ref{contrad}). This means that $y^{0}=0$,
consequently, every convergent subsequence of
$(f^{k}(x))_{k\in\mathbb{N}}$ converges to $0$. This provides that
the sequence $(f^{k}(x))_{k\in\mathbb{N}}$ is convergent to $0$,
and $x\in D_{a}(0)$.

Therefore, the ball $B(R)$ is contained in the domain of
attraction of $D_{a}(0)$.  \end{proof}

For $p\geq 0$ and $c>0$ let be $N_{p}^{c}$ the set
\begin{equation}\label{Npc.def}
    N_{p}^{c}=\{x\in\Omega: V_{p}(x)<c\}
\end{equation}
If $c=+\infty$, then $N_{p}^{c}=\Omega$.

\begin{theorem}
Let be $p\geq 0$. For any $c\in (0,(p+1)R^2]$, the set $N_{p}^{c}$
is included in the domain of attraction $D_{a}(0)$.
\label{theorem.Npc.parte.DA}
\end{theorem}

\begin{proof}
Let be $c\in (0,(p+1)R^2]$ and $x\in N_{p}^{c}$. Then
$V_{p}(x)=\sum\limits_{k=0}^{p}\|f^{k}(x)\|^{2}< c\leq (p+1)R^2$,
therefore, there exists $k\in\{0,1,..,p\}$ such that
$\|f^{k}(x)\|^{2}<R^2$. It results that $f^{k}(x)\in B(R)\subset
D_{a}(0)$, therefore, $x\in D_{a}(0)$.  \end{proof}

\begin{remark}
It is obvious that for $p\geq 0$ and $0<c'<c''$ one has
$N_{p}^{c'}\subset N_{p}^{c''}$. Therefore, for a given $p\geq 0$,
the largest part of $D_{a}(0)$ which can be found by this method
is $N_{p}^{c_{p}}$, where $c_{p}=(p+1)R^2$. In the followings, we
will use the notation $N_{p}$ instead of $N_{p}^{c_{p}}$. Shortly,
$N_{p}=\{x\in\Omega: V_{p}(x)<(p+1)R^{2}\}$ is a part of
$D_{a}(0)$. Let's note that $N_{0}=B(R)$.
\end{remark}

\begin{remark}
If $R=+\infty$ (i.e. $\Omega=\mathbb{R}^{n}$ and $\|f(x)\|<\|x\|$,
for any $x\in\mathbb{R}\setminus\{0\}$), then
$N_{p}=\mathbb{R}^{n}$ for any $p\geq 0$ and
$D_{a}(0)=\mathbb{R}^{n}$.
\end{remark}

\begin{theorem}
\label{theorem.main.norma} For the sets
$(N_{p})_{p\in\mathbb{N}}$, the following properties hold:
\begin{itemize}
    \item[a.] For any $p\geq 0$, one has $N_{p}\subset N_{p+1}$;
    \item[b.] For any $p\geq 0$ the set $N_{p}$ is invariant to $f$;
    \item[c.] For any $x\in D_{a}(0)$ there exists $p^{x}\geq 0$ such that
$x\in N_{p^{x}}$.
\end{itemize}
\end{theorem}

\begin{proof}
\emph{a.} Let be $p\geq 0$ and $x\in N_{p}$. Then
$V_{p}(x)=\sum\limits_{k=0}^{p}\|f^{k}(x)\|^{2}<(p+1)R^2$,
therefore, there exists $k\in\{0,1,..,p\}$ such that
$\|f^{k}(x)\|^{2}<R^2$. It results that $f^{k}(x)\in B(R)$ and
therefore $f^{m}(x)\in B(R)$, for any $m\geq k$. For $m=p+1$ we
obtain $\|f^{p+1}(x)\|<R$, hence
$V_{p+1}(x)=V_{p}(x)+\|f^{p+1}(x)\|^2<(p+1)R^{2}+R^{2}=(p+2)R^{2}$.
Therefore, $x\in N_{p+1}$.

\emph{b.} Let be $x\in N_{p}$. If $\|x\|<R$ then $\|f^{m}(x)\|<R$
for any $m\geq 0$ (by means of Lemma \ref{lem.B(R).proprietati}).
This implies that
$V_{p}(f(x))=\sum\limits_{k=0}^{p}\|f^{k}(f(x))\|^{2}=\sum\limits_{k=1}^{p+1}\|f^{k}(x)\|^{2}<(p+1)R^{2}$,
meaning that $f(x)\in N_{p}$.

Let's suppose that $\|x\|\geq R$. As $x\in N_{p}$, we have that
$V_{p}(x)=\sum\limits_{k=0}^{p}\|f^{k}(x)\|^{2}<(p+1)R^2$,
therefore, there exists $k\in\{0,1,..,p\}$ such that
$\|f^{k}(x)\|<R$. It results that $f^{k}(x)\in B(R)$ and therefore
$f^{m}(x)\in B(R)$, for any $m\geq k$. For $m=p+1$ we obtain
$\|f^{p+1}(x)\|<R$. This implies that
\begin{equation}
    V_{p}(f(x))=V_{p}(x)+\|f^{p+1}(x)\|^{2}-\|x\|^{2}<(p+1)R^{2}+R^{2}-R^{2}=(p+1)R^{2}
\end{equation}
therefore $f(x)\in N_{p}$.

\emph{c.} Suppose the contrary, i.e. there exist $x\in D_{a}(0)$
such that for any $p\geq 0$, $x\notin N_{p}$. Therefore,
$V_{p}(x)\geq (p+1)R^2$ for any $p\geq 0$. Passing to the limit
for $p\rightarrow\infty$ in this inequality, provides that
$V(x)=\infty$. This means $x\in\partial D_{a}(0)$ which
contradicts the fact that $x$ belongs to the open set $D_{a}(0)$.
In conclusion, there exists $p^{x}\geq 0$ such that $x\in
N_{p^{x}}$.
\end{proof}

For $p\geq 0$ let be $M_{p}=f^{-p}(B(R))=\{x\in\Omega: f^{p}(x)\in
B(R)\}$, obtained by the trajectory reversing method.

\begin{theorem}
\label{theorem.Mp.norma} The following properties hold:
\begin{itemize}
    \item[a.] $M_{p}\subset D_{a}(0)$ for any $p\geq 0$;
    \item[b.] For any $p\geq 0$, $M_{p}$ is invariant to $f$;
    \item[c.] $M_{p}\subset M_{p+1}$ for any $p\geq 0$;
    \item[d.] For any $x\in D_{a}(0)$ there
exists $p^{x}\geq 0$ such that $x\in M_{p^{x}}$.
\end{itemize}
\end{theorem}

\begin{proof}
\emph{a.} As $M_{p}=f^{-p}(B(R))$ and $B(R)\subset D_{a}(0)$ (see
Theorem \ref{theorem.B(R).parte.DA}) it is clear that
$M_{p}\subset D_{a}(0)$.

\emph{b.} and \emph{c.} follow easily by induction, using Lemma
\ref{lem.B(R).proprietati}.

\emph{d.} $x\in D_{a}(0)$ provides that $f^{p}(x)\rightarrow 0$ as
$p\rightarrow \infty$. Therefore, there exists
$p^{x}\in\mathbb{N}$ such that $f^{p}(x)\in B(R)$, for any $p\geq
p^{x}$. This provides that $x\in M_{p}$ for any $p\geq p^{x}$.
\end{proof}

Both sequences of sets $(M_{p})_{p\in\mathbb{N}}$ and
$(N_{p})_{p\in\mathbb{N}}$ are increasing, and are made up of
estimates of $D_{a}(0)$. From the practical point of view, it is
important to know which sequence converges more quickly. The next
theorem provides that the sequence $(M_{p})_{p\in\mathbb{N}}$
converges more quickly than $(N_{p})_{p\in\mathbb{N}}$, meaning
that for $p\geq 0$, the set $M_{p}$ is a larger estimate of
$D_{a}(0)$ then $N_{p}$.

\begin{theorem}
\label{theorem.incluziune.norma} For any $p\geq 0$ one has
$N_{p}\subset M_{p}$.
\end{theorem}

\begin{proof}
Let be $p\geq 0$ and $x\in N_{p}$. We have that
$V_{p}(x)=\sum\limits_{k=0}^{p}\|f^{k}(x)\|^{2}<(p+1)R^2$,
therefore, there exists $k\in\{0,1,..,p\}$ such that
$\|f^{k}(x)\|<R$. This implies that $f^{m}(x)\in B(R)$, for any
$m\geq k$. For $m=p$ we obtain $f^{p}(x)\in B(R)$, meaning that
$x\in M_{p}$.
\end{proof}

\section{Theoretical results when $A=\partial_0 f$ is a convergent matrix, i.e. $\rho(A)<1$}

As $\rho(A)\leq\|A\|$, the case when $A$ is a convergent
matrix \cite{Horn-Johnson} is more general then the case when $A$
is a contraction treated in the previous section.

\begin{proposition}
If $\rho(A)<1$, then there exists $\tilde{p}\in\mathbb{N}^{\star}$
and $r_{\tilde{p}}>0$ such that $B(r_{\tilde{p}})\subset\Omega$
and $\|f^{p}(x)\|<\|x\|$ for any
$p\in\{\tilde{p},\tilde{p}+1,..,2\tilde{p}-1\}$ and $x\in
B(r_{\tilde{p}})\setminus\{0\}$. \label{prop.ineg.raza}
\end{proposition}

\begin{proof}
We have that $\rho(A)<1$ if and only if
$\lim\limits_{p\rightarrow\infty}A^{p}=0$ (see
\cite{Horn-Johnson}), which provides that there exists
$\tilde{p}\in\mathbb{N}^{\star}$ such that $\|A^{p}\|<1$ for any
$p\geq \tilde{p}$. Let be $\tilde{p}\in\mathbb{N}^{\star}$ fixed
with this property.

The formula of variation of constants for any $p$ gives:
\begin{equation}\label{variatia.const}
    f^{p}(x)=A^{p}x+\sum_{k=0}^{p-1}A^{p-k-1}g(f^{k}(x))\qquad
    \textrm{for all }x\in\Omega\textrm{ and }p\in\mathbb{N}^{\star}
\end{equation}
Due to the fact that for any $k\in\mathbb{N}$ we have
$\lim\limits_{x\rightarrow 0}\frac{\|g(f^{k}(x))\|}{\|x\|}=0$,
there exists $r_{\tilde{p}}>0$ such that for any
$p\in\{\tilde{p},\tilde{p}+1,..,2\tilde{p}-1\}$ the following
inequality holds:
\begin{equation}\label{B(tildeR).definitie}
    \sum_{k=0}^{p-1}\|A^{p-k-1}\|\|g(f^{k}(x))\|<(1-\|A^{p}\|)\|x\|\qquad\textrm{for }x\in
    B(r_{\tilde{p}})\setminus\{0\}
\end{equation}
Let be $x\in B(r_{\tilde{p}})\setminus\{0\}$ and
$p\in\{\tilde{p},\tilde{p}+1,..,2\tilde{p}-1\}$. Using
(\ref{variatia.const}) and (\ref{B(tildeR).definitie}) we have
\begin{eqnarray}
\nonumber \|f^{p}(x)\| &=& \|A^{p}x+\sum_{k=0}^{p-1}A^{p-k-1}g(f^{k}(x))\|\leq \\
\nonumber &\leq&\|A^{p}\|\|x\|+\sum_{k=0}^{p-1}\|A^{p-k-1}\|\|g(f^{k}(x))\|<\\
   &<&(\|A^{p}\|+1-\|A^{p}\|)\|x\|=\|x\|
\end{eqnarray}
Therefore, $\|f^{p}(x)\|<\|x\|$ for
$p\in\{\tilde{p},\tilde{p}+1,..,2\tilde{p}-1\}$ and $x\in
B(r_{\tilde{p}})\setminus\{0\}$.
\end{proof}

\begin{definition}
Let be $\tilde{p}\in\mathbb{N}^{\star}$ the smallest number such
that $\|A^{p}\|<1$ for any $p\geq \tilde{p}$ (see the proof of
Proposition \ref{prop.ineg.raza}). Let be $\tilde{R}>0$ the
largest number such that $B(\tilde{R})\subset\Omega$ and
$\|f^{p}(x)\|<\|x\|$ for
$p\in\{\tilde{p},\tilde{p}+1,..,2\tilde{p}-1\}$ and $x\in
B(\tilde{R})\setminus\{0\}$.

If for any $r>0$ we have that $B(r)\subset\Omega$ and
$\|f^{p}(x)\|<\|x\|$ for any
$p\in\{\tilde{p},\tilde{p}+1,..,2\tilde{p}-1\}$ and $x\in
B(r)\setminus\{0\}$, then $\tilde{R}=+\infty$ and
$B(\tilde{R})=\Omega=\mathbb{R}^{n}$.
\end{definition}

\begin{lemma}
\begin{itemize}
    \item[a.] For any $x\in B(\tilde{R})$ and
    $p\in\{\tilde{p},\tilde{p}+1,..,2\tilde{p}-1\}$ the sequence
    $(\|f^{kp}(x)\|)_{k\in\mathbb{N}}$ is decreasing.
    \item[b.] For any $p\geq
\tilde{p}$ and $x\in B(\tilde{R})\setminus\{0\}$,
$\|f^{p}(x)\|<\|x\|$.
    \item[c.] For any $p\geq
\tilde{p}$ and $x\in B(\tilde{R})\setminus\{0\}$, $\Delta
V_{p}(x)=V_{p}(f(x))-V_{p}(x)<0$, where $V_{p}$ is defined by
(\ref{formula.Vp}).
\end{itemize}
\label{lem.B(tildeR).proprietati}
\end{lemma}

\begin{proof}
\begin{itemize}
\item[a.] If $x=0$, then $f^{p}(0)=0$, for any $p\geq 0$.

Let be $x\in B(\tilde{R})\setminus\{0\}$. We know that
$\|f^{p}(x)\|<\|x\|$ for any
$p\in\{\tilde{p},\tilde{p}+1,..,2\tilde{p}-1\}$. It results that
$f^{p}(x)\in B(\tilde{R})$ for any
$p\in\{\tilde{p},\tilde{p}+1,..,2\tilde{p}-1\}$ . This implies
that for any $k\in\mathbb{N}^{\star}$ we have
$\|f^{kp}(x)\|<\|x\|$ and $\|f^{(k+1)p}(x)\|\leq\|f^{kp}(x)\|$,
meaning that the sequence $(\|f^{kp}(x)\|)_{k\in\mathbb{N}}$ is
decreasing.

\item[b.] Let be $x\in B(\tilde{R})\setminus\{0\}$. Inequality
$\|f^{p}(x)\|<\|x\|$ is true for any
$p\in\{\tilde{p},\tilde{p}+1,..,2\tilde{p}-1\}$.

Let be $p\geq 2\tilde{p}$. There exists $q\in\mathbb{N}^{\star}$
and $p'\in\{\tilde{p},\tilde{p}+1,..,2\tilde{p}-1\}$ such that
$p=q\tilde{p}+p'$. Using \emph{a.}, we have that $f^{p'}(x)\in
B(\tilde{R})$ and $f^{q\tilde{p}}(y)\leq\|y\|$, for any $y\in
B(\tilde{R})$, therefore
\begin{equation}
    \|f^{p}(x)\|=\|f^{q\tilde{p}}(f^{p'}(x))\|\leq\|f^{p'}(x)\|<\|x\|
\end{equation}

\item[c.] results directly from \emph{b.}
\end{itemize}
\end{proof}

\begin{corollary}
For any $p\geq \tilde{p}$, there exists a maximal domain
$G_{p}\subset\Omega$ such that $0\in G_{p}$ and for any $x\in
G_{p}\setminus\{0\}$, the (positive definite) function $V_{p}$
verifies $\Delta V_{p}(x)<0$. In other words, for any $p\geq
\tilde{p}$ the function $V_{p}$ is a Lyapunov function for
(\ref{dyn.sys}) on $G_{p}$. More, $B(\tilde{R})\subset G_{p}$ for
any $p\geq \tilde{p}$.
\end{corollary}

\begin{lemma}
For any $k\geq\tilde{p}$ there exists $q_{k}\in\mathbb{N}$ such
that
\begin{equation}\label{ineg.dubla}
     \|f^{(q_{k}+3)\tilde{p}}(x)\|\leq\|f^{k}(x)\|\leq\|f^{q_{k}\tilde{p}}(x)\|\qquad\textrm{for
any }x\in B(\tilde{R})
\end{equation}
\label{lem.ineg.dubla}
\end{lemma}

\begin{proof}
Let be $k\geq\tilde{p}$. There exists a unique
$q_{k}\in\mathbb{N}$ and a unique
$p_{k}\in\{\tilde{p},\tilde{p}+1,..,2\tilde{p}-1\}$ such that
$k=q_{k}\tilde{p}+p_{k}$. Lemma \ref{lem.B(tildeR).proprietati}
provides that for any $x\in B(\tilde{R})$ we have that
$f^{q_{k}\tilde{p}}(x)\in B(\tilde{R})$ and
$\|f^{p_{k}}(x)\|\leq\|x\|$. It results that
\begin{equation}
\|f^{k}(x)\|=\|f^{p_{k}}(f^{q_{k}\tilde{p}}(x))\|\leq\|f^{q_{k}\tilde{p}}(x)\|\qquad\textrm{for
any }x\in B(\bar{R})
\end{equation}
On the other hand, we have
$(q_{k}+3)\tilde{p}=k+(3\tilde{p}-p_{k})$. As
$(3\tilde{p}-p_{k})\in \{\tilde{p}+1,\tilde{p}+2,..,2\tilde{p}\}$
and $k\geq\tilde{p}$, Lemma \ref{lem.B(tildeR).proprietati}
provides that for any $x\in B(\tilde{R})$ we have that
$f^{k}(x)\in B(\tilde{R})$ and
$\|f^{3\tilde{p}-p_{k}}(x)\|\leq\|x\|$. Therefore
\begin{equation}
    \|f^{(q_{k}+3)\tilde{p}}(x)\|=\|f^{3\tilde{p}-p_{k}}(f^{k}(x))\|\leq\|f^{k}(x)\|\qquad\textrm{for
any }x\in B(\tilde{R})
\end{equation}
Combining the two inequalities, we get that
\begin{equation}\label{ineg.dubla.folos1}
     \|f^{(q_{k}+3)\tilde{p}}(x)\|\leq\|f^{k}(x)\|\leq\|f^{q_{k}\tilde{p}}(x)\|\qquad\textrm{for
any }x\in B(\tilde{R})
\end{equation}
which concludes the proof.
\end{proof}

\begin{theorem}
    $B(\tilde{R})$ is included in the domain of attraction $D_{a}(0)$.
\label{theorem.B(tildeR).parte.DA}
\end{theorem}

\begin{proof}
Let be $x\in B(\tilde{R})\setminus\{0\}$. We have to prove that
$\lim\limits_{k\rightarrow\infty}f^{k}(x)=0$.

The sequence $(f^{k}(x))_{k\in\mathbb{N}}$ is bounded (see Lemma
\ref{lem.B(tildeR).proprietati}). Let be
$(f^{k_{j}}(x))_{j\in\mathbb{N}}$ a convergent subsequence and let
be $\lim\limits_{j\rightarrow\infty}f^{k_{j}}(x)=y^{0}$.

We suppose, without loss of generality, that $k_{j}\geq\tilde{p}$
for any $j\in\mathbb{N}$. Lemma \ref{lem.ineg.dubla} provides that
for any $j\in\mathbb{N}$ there exists $q_{j}\in\mathbb{N}$ such
that
\begin{equation}\label{ineg.dubla.folos2}
     \|f^{(q_{j}+3)\tilde{p}}(x)\|\leq\|f^{k_{j}}(x)\|\leq\|f^{q_{j}\tilde{p}}(x)\|
\end{equation}
As $(\|f^{q_{j}\tilde{p}}(x)\|)_{j\in\mathbb{N}}$ and
$(\|f^{(q_{j}+3)\tilde{p}}(x)\|)_{j\in\mathbb{N}}$ are
subsequences of the convergent sequence
$(\|f^{q\tilde{p}}(x)\|)_{q\in\mathbb{N}}$ (decreasing, according
to Lemma \ref{lem.B(tildeR).proprietati}), it results that they
are convergent. The double inequality (\ref{ineg.dubla.folos2})
provides that
$\lim\limits_{j\rightarrow\infty}\|f^{q_{j}\tilde{p}}(x)\|=\|y^{0}\|$.
Therefore,
$\lim\limits_{q\rightarrow\infty}\|f^{q\tilde{p}}(x)\|=\|y^{0}\|$.

It can be shown that
\begin{equation}\label{ineg0.raza}
    \|f^{k}(x)\|\geq \|y^{0}\| \qquad\textrm{for any } k\geq\tilde{p}
\end{equation}
For this, remark that
$\lim\limits_{q\rightarrow\infty}\|f^{q\tilde{p}}(x)\|=\|y^{0}\|$
and $(\|f^{q\tilde{p}}(x)\|)_{q\in\mathbb{N}}$ is decreasing
(Lemma \ref{lem.B(tildeR).proprietati}), which implies that
$\|f^{q\tilde{p}}(x)\|\geq \|y^{0}\|$ for any $q\in\mathbb{N}$. On
the other hand, Lemma \ref{lem.ineg.dubla} provides that for any
$k\geq\tilde{p}$ there exists $q_{k}$ such that
$\|f^{(q_{k}+3)\tilde{p}}(x)\|\leq\|f^{k}(x)\|$. Therefore,
$\|f^{k}(x)\|\geq\|f^{(q_{k}+3)\tilde{p}}(x)\|\geq\|y^{0}\|$, for
any $k\geq\tilde{p}$.

We show now that $y^{0}=0$. Suppose the contrary, i.e. $y^{0}\neq
0$.

Inequality (\ref{ineg0.raza}) becomes
\begin{equation}\label{contrad.raza}
     \|f^{k}(x)\|\geq \|y^{0}\|>0 \qquad\textrm{for any } k\geq\tilde{p}
\end{equation}
By means of Lemma \ref{lem.B(tildeR).proprietati}, we have that
$\|f^{\tilde{p}}(y^{0})\|<\|y^{0}\|$.

There exists a neighborhood $U_{f^{\tilde{p}}(y^{0})}\subset
B(\tilde{R})$ of $f^{\tilde{p}}(y^{0})$ such that for any $z\in
U_{f^{\tilde{p}}(y^{0})}$ we have $\|z\|<\|y^{0}\|$. On the other
hand, for the neighborhood $U_{f^{\tilde{p}}(y^{0})}$ there exists
a neighborhood $U_{y^{0}}\subset B(\tilde{R})$ of $y^{0}$ such
that for any $y\in U_{y^{0}}$, we have $f^{\tilde{p}}(y)\in
U_{f^{\tilde{p}}(y^{0})}$. Therefore:
\begin{equation}\label{ineg1.raza}
    \|f^{\tilde{p}}(y)\|<\|y^{0}\|\qquad \textrm{for any } y\in U_{y^{0}}
\end{equation}
As $f^{k_{j}}(x)\rightarrow y^{0}$, there exists $\bar{j}$ such
that $f^{k_{j}}(x)\in U_{y^{0}}$, for any $j\geq \bar{j}$. Making
$y=f^{k_{j}}(x)$ in (\ref{ineg1.raza}), it results that
\begin{equation}\label{ineg2.raza}
    \|f^{k_{j}+\tilde{p}}(x)\|=\|f^{\tilde{p}}(f^{k_{j}}(x))\|<\|y^{0}\|
    \qquad \textrm{for }j\geq \bar{j}
\end{equation}
which contradicts (\ref{contrad.raza}). This means that $y^{0}=0$,
consequently, every convergent subsequence of
$(f^{k}(x))_{k\in\mathbb{N}}$ converges to $0$. This provides that
the sequence $(f^{k}(x))_{k\in\mathbb{N}}$ is convergent to $0$,
and $x\in D_{a}(0)$.

Therefore, the ball $B(\tilde{R})$ is contained in the domain of
attraction of $D_{a}(0)$.  \end{proof}

\begin{theorem}
Let be $p\geq \tilde{p}$. For any $c\in (0,(p+1)\tilde{R}^2]$, the
set $N_{p}^{c}$ is included in the domain of attraction
$D_{a}(0)$. \label{theorem.Npc.parte.DA.raza}
\end{theorem}

\begin{proof}
Let be $c\in (0,(p+1)\tilde{R}^2]$ and $x\in N_{p}^{c}$. Then
$V_{p}(x)=\sum\limits_{k=0}^{p}\|f^{k}(x)\|^{2}< c\leq
(p+1)\tilde{R}^2$, therefore, there exists $k\in\{0,1,..,p\}$ such
that $\|f^{k}(x)\|^{2}<\tilde{R}^2$. It results that $f^{k}(x)\in
B(\tilde{R})\subset D_{a}(0)$, therefore, $x\in D_{a}(0)$.
\end{proof}

\begin{remark}
It is obvious that for $p\geq \tilde{p}$ and $0<c'<c''$ one has
$N_{p}^{c'}\subset N_{p}^{c''}$. Therefore, for a given $p\geq
\tilde{p}$, the largest part of $D_{a}(0)$ which can be found by
this method is $N_{p}^{\tilde{c}_{p}}$, where
$\tilde{c}_{p}=(p+1)\tilde{R}^2$. In the followings, we will use
the notation $\tilde{N}_{p}$ instead of $N_{p}^{\tilde{c}_{p}}$.
Shortly, $\tilde{N}_{p}=\{x\in\Omega:
V_{p}(x)<(p+1)\tilde{R}^{2}\}$ is a part of $D_{a}(0)$.
\end{remark}

\begin{remark}
If $\tilde{R}=+\infty$ (i.e. $\Omega=\mathbb{R}^{n}$ and
$\|f^{p}(x)\|<\|x\|$, for any
$p\in\{\tilde{p},\tilde{p}+1,..,2\tilde{p}-1\}$ and
$x\in\mathbb{R}\setminus\{0\}$), then
$\tilde{N}_{p}=\mathbb{R}^{n}$ for any $p\geq \tilde{p}$ and
$D_{a}(0)=\mathbb{R}^{n}$.
\end{remark}

\begin{theorem}
For any $x\in D_{a}(0)$ there exists $p^{x}\geq \tilde{p}$ such
that $x\in \tilde{N}_{p^{x}}$. \label{theorem.main.raza}
\end{theorem}

\begin{proof}
Let be $x\in D_{a}(0)$. Suppose the contrary, i.e. $x\notin
\tilde{N}_{p}$ for any $p\geq\tilde{p}$. Therefore, $V_{p}(x)\geq
(p+1)\tilde{R}^2$ for any $p\geq \tilde{p}$. Passing to the limit
when $p\rightarrow\infty$ in this inequality provides that
$V(x)=\infty$. This means $x\in\partial D_{a}(0)$ which
contradicts the fact that $x$ belongs to the open set $D_{a}(0)$.
In conclusion, there exists $p^{x}\geq 0$ such that $x\in
\tilde{N}_{p^{x}}$.  \end{proof}

\begin{remark}
In the case when $\|\partial_0 f\|<1$ we have shown that the
sequence of sets $(N_{p})_{p\in\mathbb{N}}$ is increasing (see
Theorem \ref{theorem.main.norma}).

\textbf{Open Question:} If $\|\partial_0 f\|\geq1$ and
$\rho(\partial_0 f)<1$, is the sequence of sets
$(\tilde{N}_{p})_{p\in\mathbb{N}}$ increasing?
\end{remark}

For $p\geq 0$ let be
$\tilde{M}_{p}=f^{-p}(B(\tilde{R}))=\{x\in\Omega: f^{p}(x)\in
B(\tilde{R})\}$, obtained by the trajectory reversing method.

\begin{theorem}
\label{theorem.Mp.raza} For the sets
$(\tilde{M}_{p})_{p\in\mathbb{N}}$ the following properties hold:
\begin{itemize}
    \item[a.] $\tilde{M}_{p}\subset D_{a}(0)$, for any $p\geq
    \tilde{p}$;
    \item[b.] $\tilde{M}_{kp}\subset \tilde{M}_{(k+1)p}$ for any $k\in\mathbb{N}$ and $p\in\{\tilde{p},\tilde{p}+1,..,2\tilde{p}-1\}$;
    \item[c.] For any $x\in D_{a}(0)$ there
exists $p^{x}\geq \tilde{p}$ such that $x\in \tilde{M}_{p^{x}}$.
\end{itemize}
\end{theorem}

\begin{proof}
\emph{a.} As $\tilde{M}_{p}=f^{-p}(B(\tilde{R}))$ and
$B(\tilde{R})\subset D_{a}(0)$ (see Theorem
\ref{theorem.B(tildeR).parte.DA}) it is clear that
$\tilde{M}_{p}\subset D_{a}(0)$.

\emph{b.} follows easily by induction, using Lemma
\ref{lem.B(tildeR).proprietati}.

\emph{c.} $x\in D_{a}(0)$ provides that $f^{p}(x)\rightarrow 0$ as
$p\rightarrow \infty$. Therefore, there exists
$p^{x}\geq\tilde{p}$ such that $f^{p}(x)\in B(\tilde{R})$, for any
$p\geq p^{x}$. This provides that $x\in \tilde{M}_{p}$ for any
$p\geq p^{x}$.
\end{proof}

\begin{remark}
The sequence of sets $(\tilde{M}_{p})_{p\in\mathbb{N}}$ is
generally not increasing (see Section 4: Numerical examples, the
Van der Pol equation).
\end{remark}

Both sequences of sets $(\tilde{M}_{p})_{p\in\mathbb{N}}$ and
$(\tilde{N}_{p})_{p\in\mathbb{N}}$ are made up of estimates of
$D_{a}(0)$. From the practical point of view, it would be
important to know which one of the sets $\tilde{M}_{p}$ or
$\tilde{N}_{p}$ is a larger estimate of $D_{a}(0)$ for a fixed
$p\geq\tilde{p}$. Such result could not be established, but the
following theorem holds:

\begin{theorem}
For any $p\geq 0$ one has $\tilde{N}_{p}\subset
\tilde{M}_{p+\tilde{p}}$.
\end{theorem}

\begin{proof}
Let be $p\geq 0$ and $x\in \tilde{N}_{p}$. We have that
$V_{p}(x)=\sum\limits_{k=0}^{p}\|f^{k}(x)\|^{2}<(p+1)\tilde{R}^2$,
therefore, there exists $k\in\{0,1,..,p\}$ such that
$\|f^{k}(x)\|<\tilde{R}$. This implies that $f^{k+m}(x)\in
B(\tilde{R})$, for any $m\geq\tilde{p}$. For $m=p-k+\tilde{p}$ we
obtain $f^{p+\tilde{p}}(x)\in B(\tilde{R})$, meaning that $x\in
\tilde{M}_{p+\tilde{p}}$.
\end{proof}

\section{Numerical examples}

\subsection{Example with known domain of attraction}
Let be the following discrete dynamical system
\begin{equation}\label{ex.known.DA}
\begin{array}{ll}
\left\{\begin{array}{l}
x_{k+1}=\frac{1}{2}x_{k}(1+x_{k}^{2}+2y_{k}^{2})\\
y_{k+1}=\frac{1}{2}y_{k}(1+x_{k}^{2}+2y_{k}^{2})
\end{array}\right.\qquad k\in\mathbb{N}
\end{array}
\end{equation}
There exists an infinity of steady states for this system: $(0,0)$
(asymptotically stable) and all the points $(x,y)$ belonging to
the ellipsis $x^{2}+2y^{2}=1$ (all unstable). The domain of
attraction of $(0,0)$ is
$D_{a}(0,0)=\{(x,y)\in\mathbb{R}^{2}:x^{2}+2y^{2}<1\}$.

As $\|\partial_{(0,0)}f\|=\frac{1}{2}$, we compute the largest
number $R>0$ such that $\|f(x)\|<\|x\|$ for any $x\in
B(R)\setminus\{0\}$, and we find $R=0.7071$.

\begin{figure}[htbp]
\includegraphics*[bb=3cm 0cm 13.5cm
10.5cm,width=12cm,angle=0]{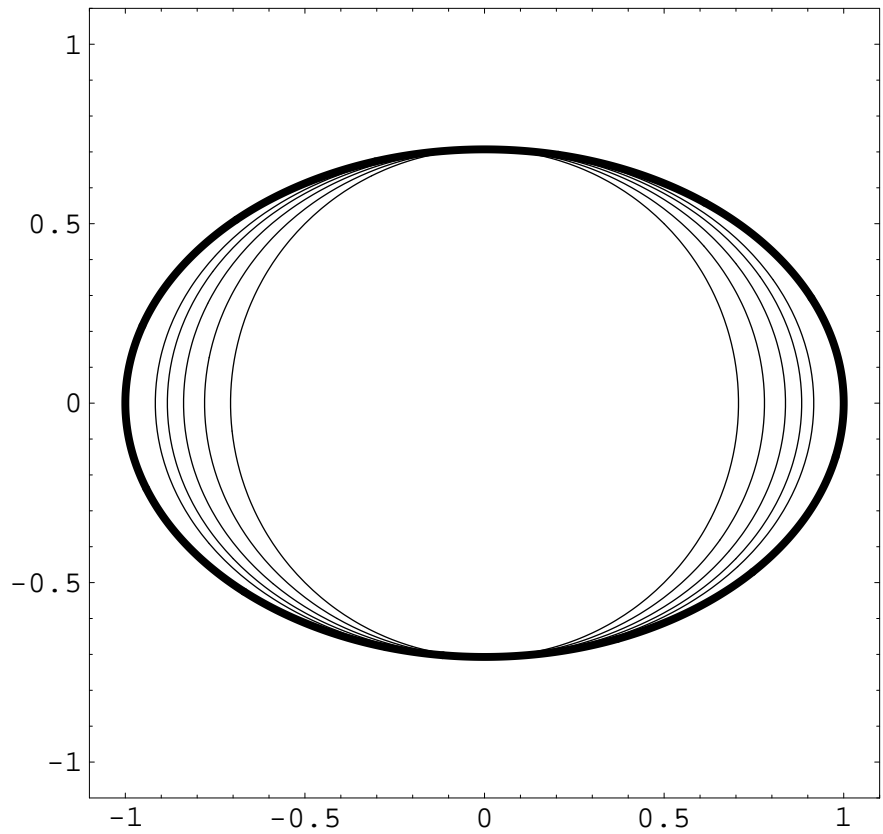} \caption{The sets $N_{p}$,
$p=\overline{0,4}$ and $\partial D_{a}(0,0)$ for
(\ref{ex.known.DA})} \label{fig.knownDA.nivel}
\end{figure}

For $p=0,1,2,3,4$ we find the $N_{p}$ sets shown in Figure
\ref{fig.knownDA.nivel}, parts of $D_{a}(0,0)$ ($N_{p}\subset
N_{p+1}$, for $p\geq 0$). In Figure \ref{fig.knownDA.nivel}, the
thick-contoured ellipsis represents the boundary of $D_{a}(0,0)$.

\begin{figure}[htbp]
\includegraphics*[bb=3cm 0cm 13.5cm
10.5cm,width=12cm,angle=0]{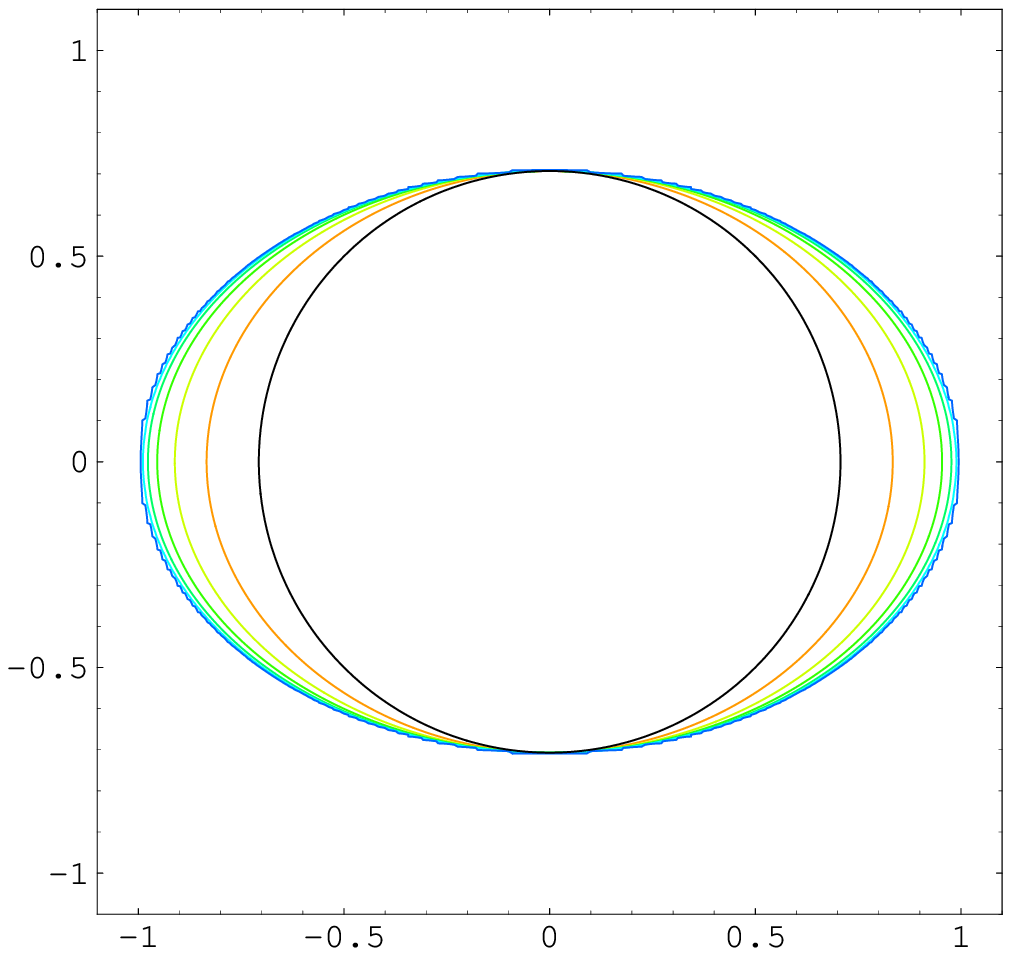}\caption{The sets
$M_{p}$, $p=\overline{0,6}$  for (\ref{ex.known.DA})}
\label{fig.knownDA.invers}
\end{figure}

In Figure \ref{fig.knownDA.invers}, the sets $M_{p}$ are
represented, for $p=\overline{0,6}$ ($M_{p}\subset M_{p+1}$, for
$p\geq 0$). Note that $M_{6}$ approximates with a good accuracy
the domain of attraction.

\subsection{Discrete predator-prey system}
We consider the discrete predator-prey system:
\begin{equation}\label{ex.predator.prey}
\begin{array}{ll}
\left\{\begin{array}{l}
x_{k+1}=ax_{k}(1-x_{k})-x_{k}y_{k}\\
y_{k+1}=\frac{1}{b}x_{k}y_{k}
\end{array}\right.
\end{array}\qquad \textrm{with }a=\frac{1}{2},b=1,k\in\mathbb{N}
\end{equation}
The steady states of this system are: $(0,0)$ (asymptotically
stable), $(-1,0)$ and $(1,-1)$ (both unstable).

We have that $\|\partial_{(0,0)}f\|=\frac{1}{2}$, and the largest
number $R>0$ such that $\|f(x)\|<\|x\|$ for any $x\in
B(R)\setminus\{0\}$ is $R=0.65$.

\begin{figure}[htbp]
\includegraphics*[bb=3cm 0cm 13.5cm
10.5cm,width=12cm,angle=0]{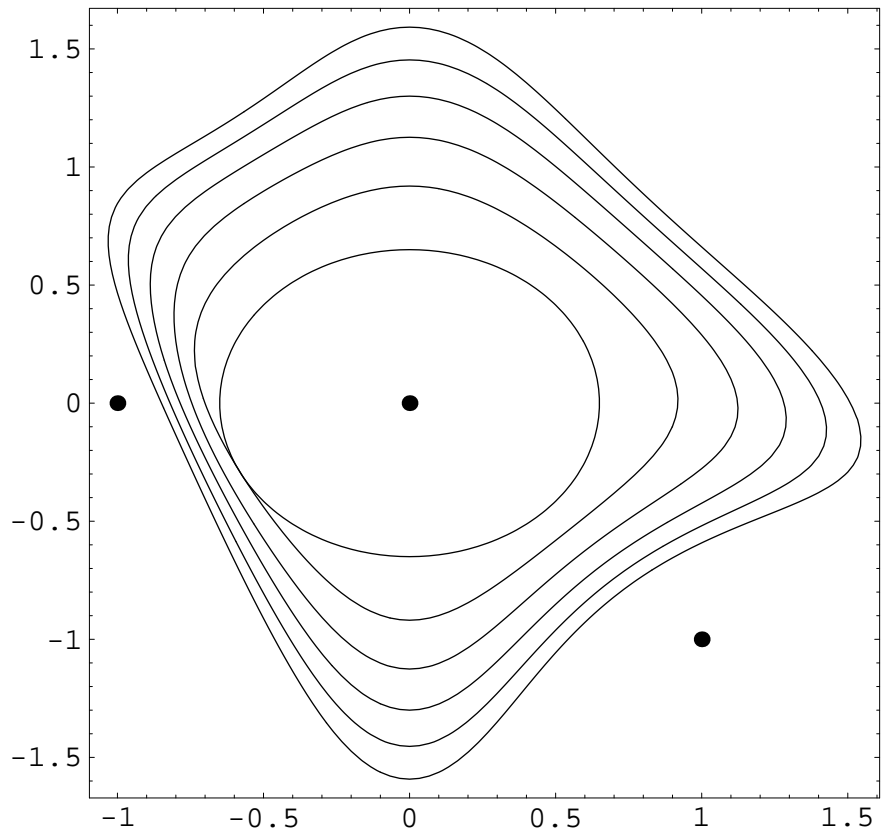}\caption{The sets
$N_{p}$, $p=\overline{0,5}$  for (\ref{ex.predator.prey})}
\label{fig.PredPrey.nivel}
\end{figure}

Figure \ref{fig.PredPrey.nivel} presents the $N_{p}$ sets for
$p=0,1,2,3,4,5$, parts of $D_{a}(0,0)$ ($N_{p}\subset N_{p+1}$,
for $p\geq 0$). The black points in Figure
\ref{fig.PredPrey.nivel} represent the steady states of the
system.

\begin{figure}[htbp]
\includegraphics*[bb=3cm 0cm 13.5cm
10.5cm,width=12cm,angle=0]{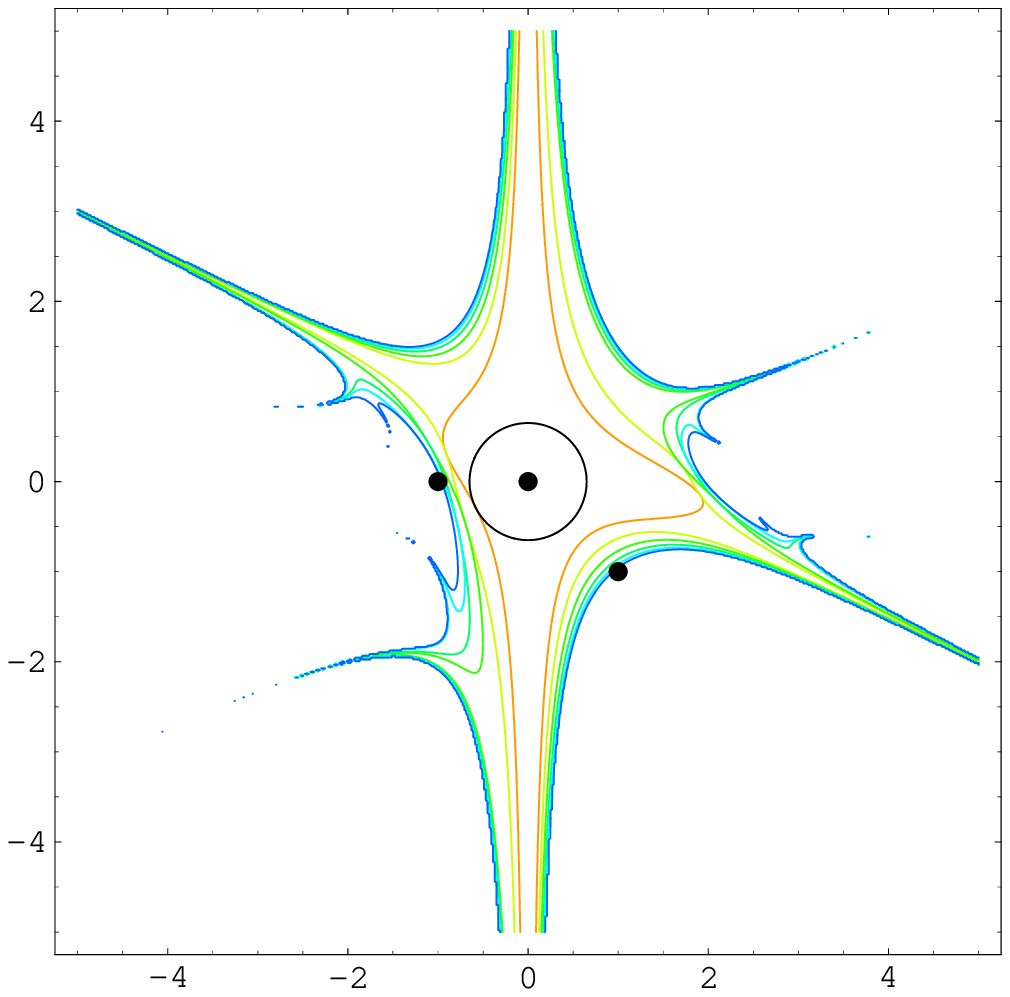}\caption{The sets
$M_{p}$, $p=\overline{0,6}$  for (\ref{ex.predator.prey})}
\label{fig.PredPrey.invers}
\end{figure}

In Figure \ref{fig.PredPrey.invers}, the sets $M_{p}$ are
represented, for $p=\overline{0,6}$ ($M_{p}\subset M_{p+1}$, for
$p\geq 0$). Note that the boundary of $M_{6}$ approaches very much
the fixed points $(-1,0)$ and $(1,-1)$, which suggests that
$M_{6}$ is a good approximation of $D_{a}(0)$.

\subsection{Discrete Van der Pol system}
\begin{figure}[htbp]
\includegraphics*[bb=3cm 0cm 13.5cm
10.5cm,width=12cm,angle=0]{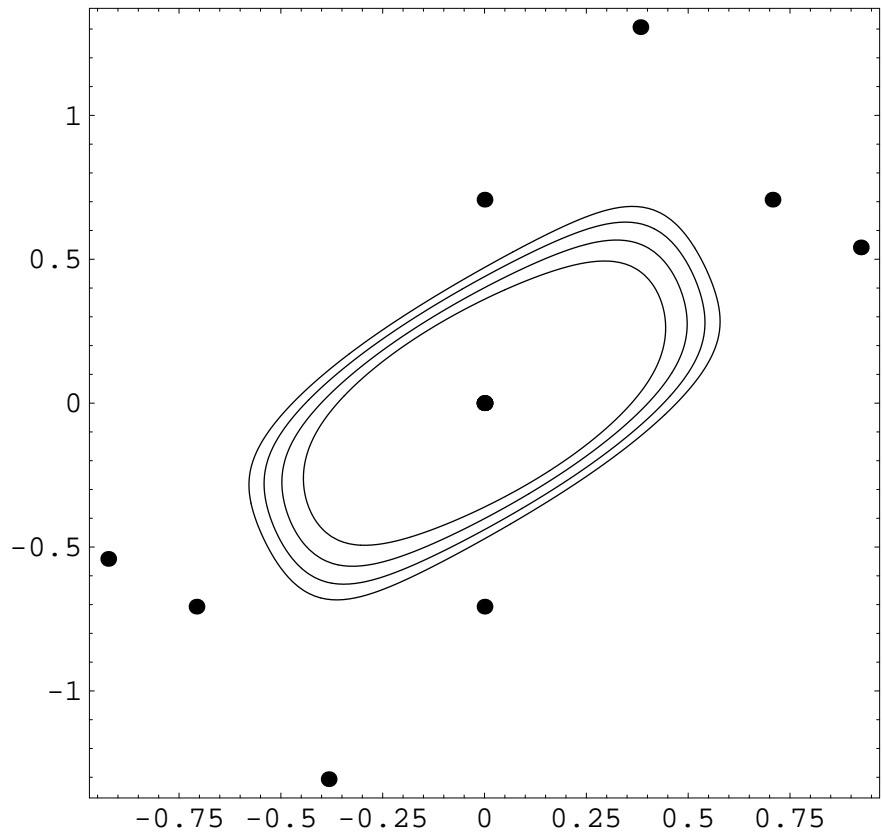}\caption{The sets
$\tilde{N}_{p}$, $p=\overline{2,5}$ for (\ref{ex.VDP})}
\label{fig.VDP.nivel}
\end{figure}

Let be the following discrete dynamical system, obtained from the
continuous Van der Pol system:
\begin{equation}\label{ex.VDP}
\begin{array}{ll}
\left\{\begin{array}{l}
x_{k+1}=x_{k}-y_{k}\\
y_{k+1}=x_{k}+(1-a)y_{k}+ax_{k}^{2}y_{k}
\end{array}\right.
\end{array}\qquad \textrm{with }a=2,k\in\mathbb{N}
\end{equation}
The only steady state of this system is $(0,0)$ which is
asymptotically stable. There are many periodic points for this
system, the periodic points of order $\overline{2,5}$ being
represented in Figure \ref{fig.VDP.nivel} by the black points.

We have that $\|\partial_{(0,0)}f\|=2$ but
$\rho(\partial_{(0,0)}f)=0$. First, we observe that for
$\tilde{p}=2$ we have that
$(\partial_{(0,0)}f)^{\tilde{p}}=O_{2}$, therefore,
$\|(\partial_{(0,0)}f)^{p}\|=0$ for any $p\geq \tilde{p}$.

The largest number $\tilde{R}>0$ such that $\|f^{p}(x)\|<\|x\|$
for $p\in\{\tilde{p},\tilde{p}+1,..,2\tilde{p}-1\}=\{2,3\}$ and
$x\in B(\tilde{R})\setminus\{0\}$ is $\tilde{R}=0.365$.

For $p=2,3,4,5$, the connected components which contain $(0,0)$ of
the $\tilde{N}_{p}$ sets are shown in Figure \ref{fig.VDP.nivel}.
We have that
$\tilde{N}_{2}\subset\tilde{N}_{3}\subset\tilde{N}_{4}\subset\tilde{N}_{5}$.

\begin{figure}[htbp]
\includegraphics*[bb=3cm 0cm 13.5cm
10.5cm,width=12cm,angle=0]{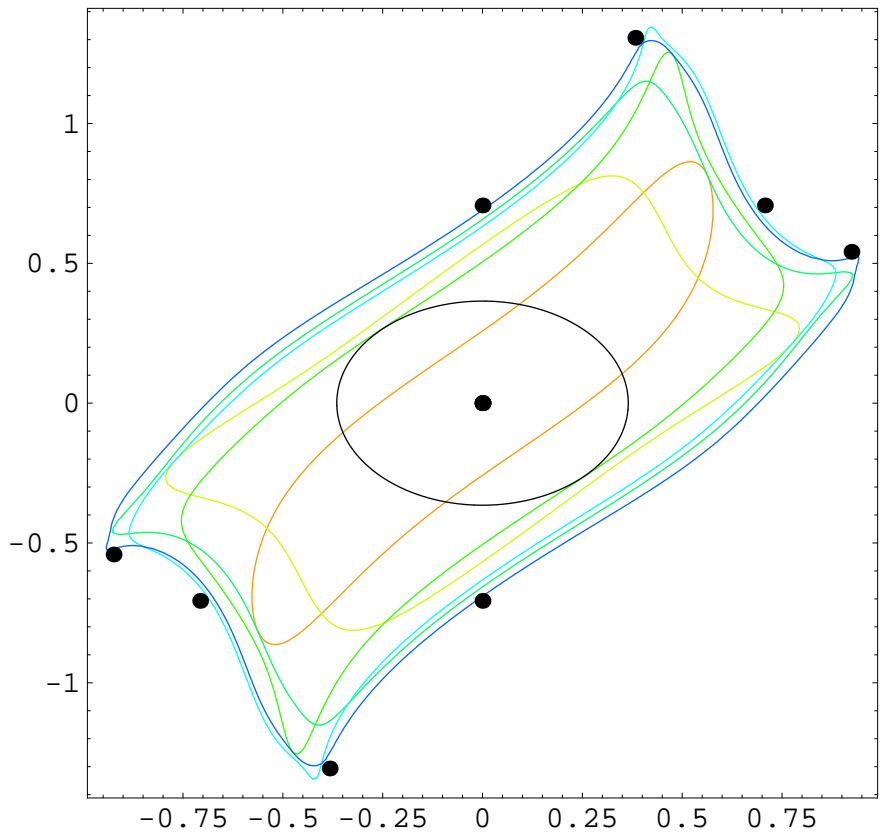}\caption{The sets
$\tilde{M}_{p}$, $p=\overline{0,6}$ for (\ref{ex.VDP})}
\label{fig.VDP.invers}
\end{figure}

In Figure \ref{fig.VDP.invers}, the sets $\tilde{M}_{p}$ are
represented, for $p=\overline{0,6}$. Note that the inclusion
$\tilde{M}_{p}\subset \tilde{M}_{p+1}$ does not hold.

\bibliography{Evabib}

\begin{thebibliography}{1}
\expandafter\ifx\csname url\endcsname\relax
  \def\url#1{\texttt{#1}}\fi
\expandafter\ifx\csname urlprefix\endcsname\relax\def\urlprefix{URL }\fi

\bibitem{Kelley-Peterson}
W.~Kelley, A.~Peterson, Differene equations, Academic Press, 2001.

\bibitem{KBBB1}
E.~Kaslik, A.~Balint, S.~Birauas, S.~Balint, Approximation of the domain of
  attraction of an asymptotically stable fixed point of a first order
  analytical system of difference equations, Nonlinear Studies 10(2)~(2) (2003)
  1--12.

\bibitem{KBGB5}
E.~Kaslik, A.~Balint, A.~Grigis, S.~Balint, An extension of the
  characterization of the domain of attraction of an asymptotically stable
  fixed point in the case of a nonlinear discrete dynamical system, in:
  Proceedings of ICNPAA 2004,
  http://www.math.uvt.ro/eng/pubs/preprints/sdea/2004/sdea2004.html.

\bibitem{Horn-Johnson}
R.~Horn, C.~Johnson, Matrix analysis, Cambridge University Press, 1985.

\end{thebibliography}

\newpage
\newpage

\newpage
\newpage
\newpage
\newpage

\end{document}